\documentclass[final,preprint,3p,sort&compress]{elsarticle}
\journal{{\tt arXiv.org}}
\usepackage[dvips]{epsfig}
\usepackage{graphicx} 
\usepackage{pdfpages}
\usepackage{latexsym}
\usepackage{verbatim}
\usepackage{amsmath,amsthm}
\usepackage{amssymb}
\usepackage{bbm}
\usepackage{fonts}
\usepackage{enumerate}
\usepackage{placeins}
\usepackage{standalone}
\usepackage{paralist}
\usepackage{subcaption}

\definecolor{greenyellow}   {cmyk}{0.15, 0   , 0.69, 0   }
\definecolor{yellow}        {cmyk}{0   , 0   , 1   , 0   }
\definecolor{goldenrod}     {cmyk}{0   , 0.10, 0.84, 0   }
\definecolor{dandelion}     {cmyk}{0   , 0.29, 0.84, 0   }
\definecolor{apricot}       {cmyk}{0   , 0.32, 0.52, 0   }
\definecolor{peach}         {cmyk}{0   , 0.50, 0.70, 0   }
\definecolor{melon}         {cmyk}{0   , 0.46, 0.50, 0   }
\definecolor{yelloworange}  {cmyk}{0   , 0.42, 1   , 0   }
\definecolor{orange}        {cmyk}{0   , 0.61, 0.87, 0   }
\definecolor{burntorange}   {cmyk}{0   , 0.51, 1   , 0   }
\definecolor{bittersweet}   {cmyk}{0   , 0.75, 1   , 0.24}
\definecolor{redorange}     {cmyk}{0   , 0.77, 0.87, 0   }
\definecolor{mahogany}      {cmyk}{0   , 0.85, 0.87, 0.35}
\definecolor{maroon}        {cmyk}{0   , 0.87, 0.68, 0.32}
\definecolor{brickred}      {cmyk}{0   , 0.89, 0.94, 0.28}
\definecolor{red}           {cmyk}{0   , 1   , 1   , 0   }
\definecolor{orangered}     {cmyk}{0   , 1   , 0.50, 0   }
\definecolor{rubinered}     {cmyk}{0   , 1   , 0.13, 0   }
\definecolor{wildstrawberry}{cmyk}{0   , 0.96, 0.39, 0   }
\definecolor{salmon}        {cmyk}{0   , 0.53, 0.38, 0   }
\definecolor{carnationpink} {cmyk}{0   , 0.63, 0   , 0   }
\definecolor{magenta}       {cmyk}{0   , 1   , 0   , 0   }
\definecolor{violetred}     {cmyk}{0   , 0.81, 0   , 0   }
\definecolor{rhodamine}     {cmyk}{0   , 0.82, 0   , 0   }
\definecolor{mulberry}      {cmyk}{0.34, 0.90, 0   , 0.02}
\definecolor{redviolet}     {cmyk}{0.07, 0.90, 0   , 0.34}
\definecolor{fuchsia}       {cmyk}{0.47, 0.91, 0   , 0.08}
\definecolor{lavender}      {cmyk}{0   , 0.48, 0   , 0   }
\definecolor{thistle}       {cmyk}{0.12, 0.59, 0   , 0   }
\definecolor{orchid}        {cmyk}{0.32, 0.64, 0   , 0   }
\definecolor{darkorchid}    {cmyk}{0.40, 0.80, 0.20, 0   }
\definecolor{purple}        {cmyk}{0.45, 0.86, 0   , 0   }
\definecolor{plum}          {cmyk}{0.50, 1   , 0   , 0   }
\definecolor{violet}        {cmyk}{0.79, 0.88, 0   , 0   }
\definecolor{royalpurple}   {cmyk}{0.75, 0.90, 0   , 0   }
\definecolor{blueviolet}    {cmyk}{0.86, 0.91, 0   , 0.04}
\definecolor{periwinkle}    {cmyk}{0.57, 0.55, 0   , 0   }
\definecolor{cadetblue}     {cmyk}{0.62, 0.57, 0.23, 0   }
\definecolor{cornflowerblue}{cmyk}{0.65, 0.13, 0   , 0   }
\definecolor{midnightblue}  {cmyk}{0.98, 0.13, 0   , 0.43}
\definecolor{navyblue}      {cmyk}{0.94, 0.54, 0   , 0   }
\definecolor{royalblue}     {cmyk}{1   , 0.50, 0   , 0   }
\definecolor{blue}          {cmyk}{1   , 1   , 0   , 0   }
\definecolor{cerulean}      {cmyk}{0.94, 0.11, 0   , 0   }
\definecolor{cyan}          {cmyk}{1   , 0   , 0   , 0   }
\definecolor{processblue}   {cmyk}{0.96, 0   , 0   , 0   }
\definecolor{skyblue}       {cmyk}{0.62, 0   , 0.12, 0   }
\definecolor{turquoise}     {cmyk}{0.85, 0   , 0.20, 0   }
\definecolor{tealblue}      {cmyk}{0.86, 0   , 0.34, 0.02}
\definecolor{aquamarine}    {cmyk}{0.82, 0   , 0.30, 0   }
\definecolor{bluegreen}     {cmyk}{0.85, 0   , 0.33, 0   }
\definecolor{emerald}       {cmyk}{1   , 0   , 0.50, 0   }
\definecolor{junglegreen}   {cmyk}{0.99, 0   , 0.52, 0   }
\definecolor{seagreen}      {cmyk}{0.69, 0   , 0.50, 0   }
\definecolor{green}         {cmyk}{1   , 0   , 1   , 0   }
\definecolor{forestgreen}   {cmyk}{0.91, 0   , 0.88, 0.12}
\definecolor{pinegreen}     {cmyk}{0.92, 0   , 0.59, 0.25}
\definecolor{limegreen}     {cmyk}{0.50, 0   , 1   , 0   }
\definecolor{yellowgreen}   {cmyk}{0.44, 0   , 0.74, 0   }
\definecolor{springgreen}   {cmyk}{0.26, 0   , 0.76, 0   }
\definecolor{olivegreen}    {cmyk}{0.64, 0   , 0.95, 0.40}
\definecolor{rawsienna}     {cmyk}{0   , 0.72, 1   , 0.45}
\definecolor{sepia}         {cmyk}{0   , 0.83, 1   , 0.70}
\definecolor{brown}         {cmyk}{0   , 0.81, 1   , 0.60}
\definecolor{tan}           {cmyk}{0.14, 0.42, 0.56, 0   }
\definecolor{gray}          {cmyk}{0   , 0   , 0   , 0.50}
\definecolor{black}         {cmyk}{0   , 0   , 0   , 1   }
\definecolor{white}         {cmyk}{0   , 0   , 0   , 0   } 

 \usepackage[]{hyperref} 
        
\newcommand{\externaltikz}[2]{\includegraphics{Externals/#1}}
 \usepackage{relinput}

\newcounter{tikzsubfigcounter}[figure]
\renewcommand{\thetikzsubfigcounter}{\thesection.\the\numexpr\value{figure}+1\relax\alph{tikzsubfigcounter}}

\newcounter{tikzsubfigcounterinvisible}[figure]
\renewcommand{\thetikzsubfigcounterinvisible}{\thesection.\the\numexpr\value{figure}+1\relax\alph{tikzsubfigcounterinvisible}}

\numberwithin{equation}{section}

\newcommand{\bdm}{\begin{displaymath}}
\newcommand{\edm}{\end{displaymath}}
\newcommand{\beq}{\begin{equation}}
\newcommand{\eeq}{\end{equation}}
\newcommand{\beqa}{\begin{eqnarray}}
\newcommand{\eeqa}{\end{eqnarray}}

\title{Partial-moment minimum-entropy models for kinetic chemotaxis equations in one and two dimensions}
\author[jr]{Juliane Ritter}
\author[ak]{Axel Klar}
\author[fs]{Florian Schneider}
\address[jr]{Fachbereich Mathematik, TU Kaiserslautern, Erwin-Schr\"odinger-Str., 67663 Kaiserslautern, Germany, {\tt j\_ritter@rhrk.uni-kl.de}}
\address[ak]{Fachbereich Mathematik, TU Kaiserslautern, Erwin-Schr\"odinger-Str., 67663 Kaiserslautern, Germany, {\tt klar@mathematik.uni-kl.de}}
\address[fs]{Fachbereich Mathematik, TU Kaiserslautern, Erwin-Schr\"odinger-Str., 67663 Kaiserslautern, Germany, {\tt schneider@mathematik.uni-kl.de}}

\date{}

\newlength{\figureheight}
\newlength{\figurewidth}
\newcommand{\R}{\mathbb{R}}
\newcommand{\norm}[2]{\ensuremath{\left\Vert #1 \right\Vert_{#2}}}
\newcommand{\Lp}[1]{\ensuremath{L_{#1}}}
\newcommand{\basis}[1][ ]{{\ensuremath{\bb_{#1}}}} 
\newcommand{\moments}[1][ ]{\ensuremath{\bu_{#1}}} 
\newcommand{\x}{\ensuremath{\bx}} 
\newcommand{\ints}[1]{\ensuremath{\left<#1\right>}} 
\newcommand{\ansatz}[1][ ]{\ensuremath{\hat{f}_{#1}}}
\newcommand{\PN}[1][N]{\ensuremath{\text{P}_{#1}}}
\newcommand{\MN}[1][N]{\ensuremath{\text{M}_{#1}}}
\newcommand{\HMN}[1][N]{\ensuremath{\text{HM}_{#1}}}
\newcommand{\HPN}[1][N]{\ensuremath{\text{HP}_{#1}}}
\newcommand{\QPN}[1][N]{\ensuremath{\text{QP}_{#1}}}
\newcommand{\QMN}[1][N]{\ensuremath{\text{QM}_{#1}}}
\def\quand{\quad \mbox{and} \quad}
\newcommand{\entropy}{\ensuremath{\eta}} 
\newcommand{\entropyFunctional}{\ensuremath{\mathcal{H}}} 
\newcommand{\indicator}[1]{\ensuremath{\mathbbm{1}_{#1}}}

\begin{document}

\begin{abstract}
The aim of this work is to investigate the application of partial moment approximations to kinetic chemotaxis equations  in one and two spatial dimensions. Starting with  a kinetic equation for the cell  densities we apply a half-/quarter-moments method with different  closure relations to derive macroscopic equations. Appropriate  numerical schemes are presented as well as numerical results for several test cases. The resulting  solutions are compared  to kinetic reference solutions and solutions computed using a full moment method with a linear superposition strategy.
\end{abstract}
\maketitle

\noindent


\section{Introduction}
\label{sec:Introduction}

The migration of cells is a complex process that is influenced by many factors such as external light, the pH or the oxygen concentration. In the following we concentrate on chemotaxis, the movement of cells in response to a chemical stimulus. A substance which causes cells to move in the direction of its gradient is called chemoattractant. The ideas of this work can easily be extended to chemorepellants, which have the opposite effect. 
Chemotaxis plays an important role in a lot of biological processes: It drives the movement of bacteria towards food and away from poisons or leads the sperm in the direction of the egg during fertilization. In multi-cellular organisms, it controls the guided accumulation of cells during embryological development and the movement of lymphocytes in the process of immunological response~\cite{bellomo, chalub, natalini}. During cancer metastasis, mechanisms that allow chemotaxis can be subverted. Therefore a better understanding of the associated processes may lead to the development of novel therapeutic strategies. 
The movement of many bacteria, such as Escherichia coli which has been studied and described most intensely, is controlled by the alignment of their flagella, whip-shaped appendices with a rotary motor at their bases that are embedded in the cell membrane. Counter-clockwise rotation aligns the flagella, causing the bacterium to swim in a straight line (``run phase''); clockwise rotation causes the flagella to point in different directions, resulting in a movement on the spot (``tumble phase''). The latter re-orients the bacterium, so that overall we observe a random walk~\cite{alt, tindall}. A chemical stimulus influences the motion in the following way: if receptors sense that the bacterium is moving in the direction of the chemoattractant gradient, the ``run phase'' will be extended; if the concentration of the chemoattractant is decreasing in the direction of movement, it will be shortened. This results in a biased random walk~\cite{chalub}. Consideration of additional effects on the chemoattractant, such as production by the bacteria themselves, decay and diffusion, further increases the complexity of the model. 

The original equations to model chemotaxis are the  Keller-Segel equations.
These equations and in particular, the properties of their solutions  have been intensively investigated, see for example
 \cite{CA06,HE97,KS70,RA95,BO10,KU12,HI09,KS71}.
For a survey and an extended reference list, see for example
 \cite{BBNS12}.

Our starting point is the  classical kinetic chemotaxis  equation \cite{CMPS04}.
Scaling it with the so called  diffusive scaling leads to the  Keller-Segel equation, see \cite{CMPS04}. 
In general, the derivation of  Keller-Segel type models, including flux-limited diffusion models and  Fokker-Planck type  models, from underlying 
kinetic or microscopic models is discussed for example in 
 \cite{BBNS12,Chav06,CMPS04}.
In particular, using moment closure approaches one may obtain macroscopic equations
intermediate between kinetic and Keller-Segel equations, see the above mentioned references
or \cite{Filbet,HI09}.
First order full moment equations with a linear closure function are sometimes called the   Cattaneo equations.
Applying maximum entropy closures one obtains improved first order full moment models \cite{Levermore1996}.
 Half-  and  quarter-moment closures have been developed for   kinetic radiative  transfer equations  in \cite{Frank2006,FraDubKla04,DubKla02}.
It has been shown for these applications in various numerical experiments that the partial moment moment methods yields macroscopic models that can produce satisfying approximations. We refer to \cite{Frank2006}.

The aim of this work is to investigate the application of partial moment approximations to kinetic chemotaxis equations  in one and two dimensions. Starting with  a kinetic equation for the cell  densities we apply a half- and quarter-moment method with varying closure relations to derive macroscopic equations  in section 2. By applying numerical schemes that use certain properties of the moment systems and  implementing  it using Matlab \cite{MATLAB:2012}, we obtain numerical results for several different test cases. Moreover, we compare the  results to kinetic reference solutions and solutions computed using a full moment method and a linear superposition strategy.
This work should be seen as a first step towards an efficient simulation strategy for kinetic chemotaxis equations via further refinement of the sphere, leading from a quarter moment model to a general first-order partial-moment model, which will hopefully converge to the true kinetic solution with a small number of refinements.
This strategy should yield similiar results as higher order  moment methods while being numerically much more efficient due to the inherent structur of the first order partial moment models.

\section{Chemotaxis equations}

The dynamics of chemotaxis can be modelled by the kinetic equations
\begin{align}
\partial_t f + v \cdot \nabla_\x f &= -\lambda(f - C_V \rho) + C_V \alpha \rho v \cdot \Phi(\nabla m) \label{eq:cell}\\
\partial_t m - D_m \Delta m &= \beta\cdot\rho - \delta m \label{eq:chem}
\end{align}
where $f(t,\x,v)$ is the density of cells at time $t$ and location $\x\in\R^d$, with velocity $v\in V$, whereas $\rho(t,\x) = \int_V f~dv$ describes the overall density of cells at time $t$ and location $\x$. $m(t,\x)$ is the concentration of the chemoattractant at time $t$ and location $\x$. $V$ is the set of admissible velocities; since we assume that the cells move in arbitrary directions, but with constant speed, we have $V=S^2 = \{v\in\R^3~|~ \norm{v}{2}=1\}$ in three dimensions and consider the projections $V=\left[-1,1\right]$ in one and $V=B_1(0) =  \{v\in\R^2~|~ \norm{v}{2}\leq 1\}$ in two dimensions~\cite{bellomo2}. The normalization constant $C_V$ is determined by $V$: $C_V=\frac{1}{2}$ in one and $C_V=\frac{1}{4\pi}$ in two dimensions. The remaining coefficients characterize the biological system: $\lambda$ and $\alpha$ describe the diffusivity and the chemotactic sensitivity of the cells. $D_m$ is the diffusivity, $\beta$ the production rate by the cells and $\delta$ the rate of chemical decay of the chemoattractant. The function $\Phi$ acts as a limiter for the influence of the chemoattractant gradient $\nabla m$, which models the fact, that the ``run phase'' can only be extended to a certain extent. In the following we use 
\begin{equation*}
\Phi(\x) = \begin{cases}
(\frac{\|\x\|-s}{\sqrt{1+(\|\x\|-s)^2}}+s)\frac{\x}{\|\x\|} & \|\x\| \geq s \\
\x & \|\x\| \leq s
\end{cases} 
\end{equation*}
where the parameter $s \geq 0$ determines the extent of the limiting: $\max(\|\Phi(\x)\|) \leq s+1$. Using this limiter in our simulation prevents blow up of the solution in finite time~\cite{chertock}.

Assuming that $\lambda \geq C_V\alpha\left(s+1\right)$, the right-hand side of \eqref{eq:cell} can be written in the turning-kernel representation with non-negative kernel, which ensures that \eqref{eq:cell} admits a non-negative solution $f$ \cite{chalub}.

We note that equation \eqref{eq:cell} is related via a diffusive scaling limit
$t \rightarrow \epsilon^2 t$ and $x \rightarrow \epsilon x$
to a special case of the general Patlak-Keller-Segel
\begin{equation*}
\partial_t \rho - \frac{1}{3\lambda} \Delta \rho = -\frac{\alpha}{3\lambda} \nabla \cdot (\rho \Phi(\nabla m)) \ .
\end{equation*}
See \cite{hillen} or \cite{chalub} for details and rigorous proofs. 

\section{Moment models}

In this section we introduce the method of moments and explore how it can be used to derive macroscopic equations from our kinetic equation~\eqref{eq:cell}.
It can be seen as a Galerkin approximation in the velocity component $v$, by projecting the kinetic density $f(t,x,v)$ in $v$ onto a finite-dimensional subspace of $\Lp{2}(V,\R)$. Assume that this subspace is spanned by the basis $\basis:V\to \R^n$, moments of $f$ are defined as
\begin{align*}
\moments(t,\x) = \int_V \basis(v) f(t,\x,v)~dv =: \ints{\basis f},
\end{align*}
where the integration is meant componentwise. Equations for $\moments$ can be obtained by multiplying \eqref{eq:cell} with $\basis$ and integrating over $V$, giving 
\begin{align}
\partial_t \moments + \ints{v \cdot \nabla_x \basis \ansatz} &= -\lambda(\moments - C_V \rho\ints{\basis}) + C_V \alpha \rho \ints{v \cdot \Phi(\nabla m)\basis} \label{eq:moments}.
\end{align}
Since the product of the components of $v$ and $\basis$ are generally not in the span of $\basis$, the ansatz $\ansatz$ has to be specified, determining the closure relation of \eqref{eq:moments}.
A typical example for $\basis$ is the set of monomials in $v$, i.e. $\basis = \left(1,v,v\otimes v,v^{\otimes 3},\ldots v^{\otimes n}\right)$, leading to what is commonly referred to as $\PN$ and $\MN$ models, depending on the choice of $\ansatz$ \cite{Lewis-Miller-1984,Lev84,Levermore1996}.
Since full-moment models (i.e. integration over the whole velocity domain $V$) suffer from severe problems \cite{Hauck2010}, we will focus on partial moment models, which build moments for partitions of $V$ separately \cite{SchFraPin04,FraDubKla04,DubKla02}.
Since some of the derived models may contain unphysical negative densities $\rho$ we replace \eqref{eq:chem} with
\begin{align*}
\partial_t m - D_m \Delta m &= \beta\cdot\max(\rho,0) - \delta m.
\end{align*}
This ensures that $m\geq 0$ is still guaranteed.

\subsection{One dimension}
Given a density function $f(t,x,v)$ with $t \in \mathbb{R}^+$, $x \in D \subset \mathbb{R}$ and \mbox{$v \in V=\left[-1,1\right]$}, we define the zeroth, first and second half-moments as
\begin{align}
\rho_\pm (t,x) := \int_{V_\pm} f(t,x,v) dv,\quad 
q_\pm (t,x) := \int_{V_\pm} v f(t,x,v) dv,\quad
r_\pm (t,x) := \int_{V_\pm} v^2 f(t,x,v) dv,\label{eq:Moments1D}
\end{align}
with $V_-:=\left[-1,0\right]$, $V_+:=\left[0,1\right]$.
Further assuming that $f\geq 0$, the monotonicity of the integral directly implies the sign restrictions $\rho_\pm \geq 0$, $\pm q_\pm \geq 0$ and $r_\pm \geq 0$.
We can derive further properties for the normalized first and second half-moments
\begin{align*}
u_\pm := \frac{q_\pm}{\rho_\pm} \quand
w_\pm := \frac{r_\pm}{\rho_\pm}\,.
\end{align*}
Since $|v| \leq 1$ for all $v \in V$, we have
\begin{equation*}
\rho_\pm = \int_{V_\pm} f dv \geq \int_{V_\pm} |v|f dv \geq \left| \int_{V_\pm} v f dv \right| = |q_\pm|.
\end{equation*}
Therefore $|u_\pm| \leq 1$, which implies $u_\pm \in V_\pm$ by using the sign restrictions. Additionally, we see
\begin{align*}
\pm q_\pm &= \int_{V_\pm} v f dv \geq \int_{V_\pm} v^2 f dv = r_\pm\,.
\end{align*}
Using the Cauchy-Schwartz inequality we also observe
\begin{equation*}
q_\pm^2 = \left(\int_{V_\pm} v f dv \right)^2 \leq \left(\int_{V_\pm} v^2 f dv \right) \left(\int_{V_\pm} f dv \right) = r_\pm \rho_\pm.
\end{equation*}
Thus, the normalized second moments satisfy 
\begin{align}
\label{eq:Realizability1D}
u_\pm^2 \leq w_\pm\leq \pm u_\pm\,.
\end{align}
The relations derived above are necessary and sufficient for the existence of a non-negative density function $f$ realizing a set of moments $\rho_\pm$, $q_\pm$ and $r_\pm$. For that reason they are also called realizability conditions. For a proof of sufficiency, see e.g. \cite{Curto1991,Schneider2014}.
In one dimension the kinetic equation~\eqref{eq:cell} for the cell density $f$ simplifies to
\begin{equation*}
\partial_t f + v \partial_x f = -\lambda(f - \frac{1}{2} \rho) + \frac{1}{2} \alpha \rho v \phi(\partial_x m),
\end{equation*}
where we used $C_V=\left(\int_{V} 1 dv\right)^{-1} = \frac{1}{2}$. \\
Integrating over $V_\pm$ yields
\begin{align}
\partial_t \rho_\pm + \partial_x q_\pm &= \frac{ \lambda}{2}\left(\rho_++\rho_- -2\rho_\pm\right) + \frac{\alpha}{4} \left(\rho_++\rho_-\right) \Phi(\partial_x m) \label{eq:1D_PDE_first}.
\end{align}
Additionally, multiplication by $v$ and integration over $V_\pm$ gives
\begin{align}
\partial_t q_\pm + \partial_x r_\pm &= \pm\frac{ \lambda}{4} \left(\rho_++\rho_- \mp 4q_\pm\right) + \frac{\alpha}{6}  \left(\rho_++\rho_-\right) \Phi(\partial_x m) \label{eq:1D_PDE_last}.
\end{align}
Here, we used the obvious properties $\int_{V_\pm} 1 dv=1$, $\int_{V_\pm} v dv=\pm \frac{1}{2}$, $\int_{V_\pm} v^2 dv = \frac{1}{3}$.
This system of four partial differential equations~\eqref{eq:1D_PDE_first}-\eqref{eq:1D_PDE_last} in six variables $\rho_\pm$, $q_\pm$ and $r_\pm$ is under-determined, requiring to derive closure relations $r_\pm=r_\pm(\rho_\pm,q_\pm)$.

\subsubsection{Linear closure}
A very cheap closure relation can by derived by using the linear ansatz
\begin{equation*}
\ansatz(t,x,v)=
\begin{cases}
a_+(t,x) + b_+(t,x) v &\text{ if } v \in V_+, \\
a_-(t,x) + b_-(t,x) v &\text{ if } v \in V_-.
\end{cases}
\end{equation*}
Plugging this into \eqref{eq:Moments1D} yields
\begin{align*}
r_\pm &= -\frac{1}{6} \rho_\pm \mp q_\pm.
\end{align*}
This simple closure relation can be implemented very easily but contradicts the realizability conditions \eqref{eq:Realizability1D} since e.g. for $\rho_\pm = 1$ and $q_\pm = 0$ it follows that $\rho_\pm r_\pm = -\frac16<0 = q_\pm^2$. Thus, the approximation will behave physically wrong (compared to the original kinetic equation) in these regimes, where the underlying ansatz is negative. This is similar to the case of the $\PN$ equations, see e.g. \cite{Schneider2014,Klar2014,Hauck2010}. This model will be called $\QPN[1]$.

\subsubsection{Minimum-entropy closure}
\begin{figure}[htbp]
	\centering
	\externaltikz{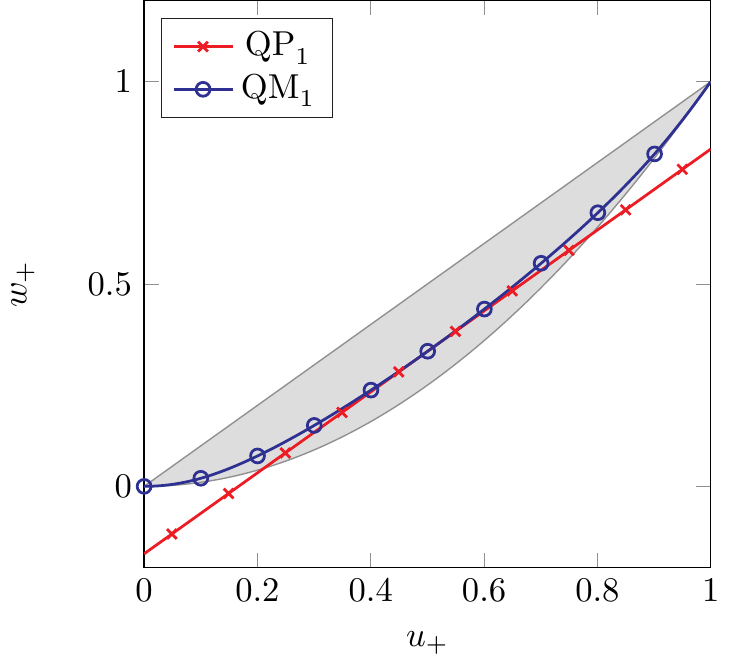}{\input{images/1D_closures}}
	\caption{Comparison of different closure relations. The set of realizable moments is shown in grey.}
	\label{fig:1D_closures}
\end{figure}

The exponential ansatz
\begin{equation*}
\ansatz(t,x,v)=
\begin{cases}
\exp(a_+(t,x) + b_+(t,x) v) &\text{ if } v \in V_+, \\
\exp(a_-(t,x) + b_-(t,x) v) &\text{ if } v \in V_-
\end{cases}
\end{equation*}
can be derived as the solution of the constrained minimization problem
\begin{align}
 \label{eq:entropyFunctional}
\min \entropyFunctional(\ansatz) = \ints{\entropy(\ansatz)}
 \end{align}
 under the moment constraints
 \begin{align}
 \label{eq:MomentConstraints}
 \ints{\basis\ansatz} = \moments,
 \end{align}
 where $\entropy$ is given by the Maxwell-Boltzmann entropy $\entropy(\ansatz) = \ansatz\log\left(\ansatz\right)-\ansatz$ \cite{Levermore1996}. The basis is chosen to be $\basis = \left(\indicator{V_-},v\indicator{V_-},\indicator{V_+},v\indicator{V_+}\right)$ with the indicator function $\indicator{}$ on $V_\pm$.
Again, calculating the moment integrals \eqref{eq:Moments1D} for $\ansatz$ yields the following relations for the normalized first and second half-moments:
\begin{alignat*}{3}
u_\pm &= \frac{q_\pm}{\rho_\pm} &&= \frac{(\pm b_\pm - 1) \exp(\pm b_\pm) + 1}{b_\pm (\exp(\pm b_\pm) - 1)}, \\
w_\pm &= \frac{r_\pm}{\rho_\pm} &&=\frac{(b_\pm^2 \mp 2b_\pm + 2) \exp(\pm b_\pm) - 2}{b_\pm^2 (\exp(\pm b_\pm) - 1)}.
\end{alignat*}
Since these expressions depend only on $b_\pm$, not on $a_\pm$, we can invert them numerically using tabulation, see e.g. \cite{Frank2005,Frank07}. It is more computationally expensive than the linear closure relation, but satisfies the realizability conditions \eqref{eq:Realizability1D}, see Figure~\ref{fig:1D_closures}. This model will be called $\QMN[1]$.

\subsection{Two dimensions}
In this section the definitions and methods of the previous section are extended to two dimensions. Most calculations will work analogously.
For a given density $f(t,\x,v)$ with $t \in \mathbb{R}^+$, $\x:=(x,y)^T \in D \subset \R^2$ and \mbox{$v:=(v_x,v_y)^T \in V=B_1(0)$}, we define quarter-moments up to second order as
\begin{align*}
\rho_{\pm \pm'}(t,\x) &:= \int_{V_{\pm \pm'}} f(t,\x,v(\phi,r)) d\phi dr =: \ints{f}_{\pm \pm'}, \\
q_{\pm \pm'}(t,\x)&:= \int_{V_{\pm \pm'}} v(\phi,r) f(t,\x,v(\phi,r)) d\phi dr =: \begin{pmatrix}
q_{\pm \pm'}^{x}(t,\x)\\q_{\pm \pm'}^{y}(t,\x) 
\end{pmatrix}, \\
r_{\pm \pm'}(t,\x) &:= \int_{V_{\pm \pm'}} v(\phi,r)\otimes v(\phi,r) f(t,\x,v(\phi,r)) d\phi dr =: \begin{pmatrix}
r_{\pm \pm'}^{xx}(t,\x) & r_{\pm \pm'}^{xy}(t,\x)\\r_{\pm \pm'}^{xy}(t,\x) & r_{\pm \pm'}^{yy}(t,\x)
\end{pmatrix},
\end{align*}
with $v(\phi,r):= (\sqrt{1-r^2} \cos(\phi),\sqrt{1-r^2} \sin(\phi))^T$ and $V_{++}:=[-1,1]\times[0,\frac{\pi}{2}]$, $V_{-+}:=[-1,1]\times[\frac{\pi}{2},\pi]$, $V_{--}:=[-1,1]\times[\pi,\frac{3\pi}{2}]$  $V_{+-}:=[-1,1]\times[\frac{3\pi}{2},2\pi]$. 
Again, normalized moments can be defined as 
\begin{align*}
u_{\pm \pm'} := \frac{q_{\pm \pm'}}{\rho_{\pm \pm'}} \quand
w_{\pm \pm'} := \frac{r_{\pm \pm'}}{\rho_{\pm \pm'}}\,,
\end{align*}
satisfying similar relations as in one dimension. For example, the first moment has to be located within the corresponding quarter sphere, i.e. $u_{\pm \pm'}\in v(V_{\pm \pm'})$ \cite{Schneider2015c,Frank2006,Ritter2015}.
Equations for the moments $\moments$ can be obtained by testing \eqref{eq:cell} with $1$ and $v$, and integrating over the four quarter spheres $V_{\pm \pm'}$:
\begin{align}
\partial_t \rho_{\pm \pm'} + \nabla\cdot q_{\pm \pm'} &= \frac{\lambda}{4} \left(\rho-4\rho_{\pm \pm'}\right) + \frac{1}{4\pi} \alpha \rho \ints{v}_{\pm \pm'}\cdot\Phi(\nabla m) \label{2D_PDE_first}, \\
\partial_t q_{\pm \pm'} + \nabla\cdot r_{\pm \pm'} &= \frac{\lambda}{8} \left(\rho-8q_{\pm \pm'}\right) + \frac{1}{4\pi} \alpha \rho \ints{v\otimes v}_{\pm \pm'} \Phi(\nabla m).\label{2D_PDE_last}
\end{align}

\begin{figure}[htbp]
	\centering
	\externaltikz{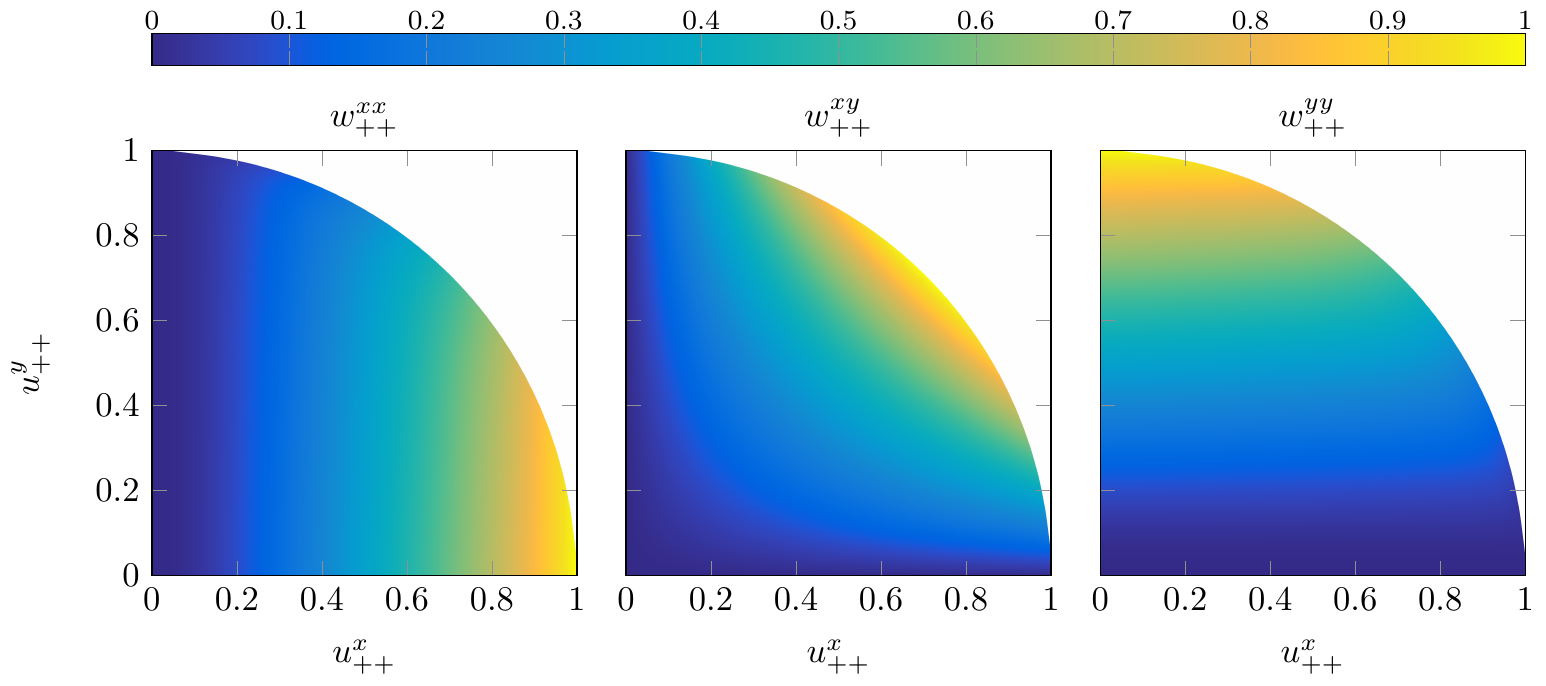}{\input{images/QM1_second_moment_cell/QM1_second_moment_cell}}
	\caption{Quarter-moment minimum-entropy closure.}
	\label{fig:QM1_second_moment_cell}
\end{figure}

The linear and minimum-entropy closures can be derived as in one dimension, using the ansätze
\begin{equation*}
\ansatz =
\begin{cases}
a_{++} + b_{++} \cdot v &\text{ if } (\phi,r) \in V_{++}, \\
a_{-+} + b_{-+} \cdot v &\text{ if } (\phi,r) \in V_{-+}, \\
a_{--} + b_{--} \cdot v &\text{ if } (\phi,r) \in V_{--}, \\
a_{+-} + b_{+-} \cdot v &\text{ if } (\phi,r) \in V_{+-}, \ .
\end{cases}
\end{equation*}
and
\begin{equation*}
\ansatz =
\begin{cases}
\exp\left(a_{++} + b_{++} \cdot v\right) &\text{ if } (\phi,r) \in V_{++}, \\
\exp\left(a_{-+} + b_{-+} \cdot v\right) &\text{ if } (\phi,r) \in V_{-+}, \\
\exp\left(a_{--} + b_{--} \cdot v\right) &\text{ if } (\phi,r) \in V_{--}, \\
\exp\left(a_{+-} + b_{+-} \cdot v\right) &\text{ if } (\phi,r) \in V_{+-},
\end{cases}
\end{equation*}
respectively. The minimum-entropy closure can be obtained again using a suitable two-dimensional table-lookup \cite{Frank2006}. Its second moment is depicted in Figure~\ref{fig:QM1_second_moment_cell}.

\section{Numerical scheme}
\begin{figure}[htbp]
	\begin{subfigure}{0.49\linewidth}
		\externaltikz{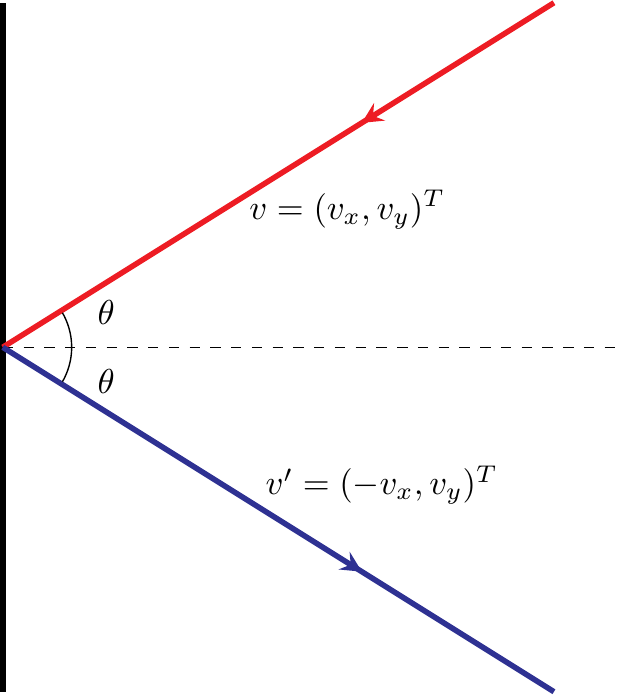}{\input{images/reflection_1}}
		\caption{$v=v(\phi,r)$ with $(\phi,r) \in V_{--}$}
	\end{subfigure}
	\begin{subfigure}{0.49\linewidth}
		\externaltikz{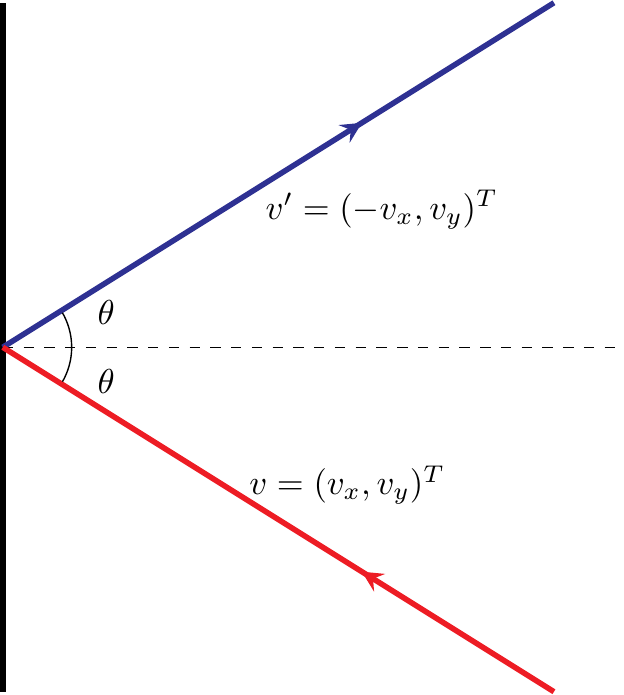}{\input{images/reflection_2}}
		\caption{$v=v(\phi,r)$ with $(\phi,r) \in V_{-+}$}
	\end{subfigure}
	\caption{Illustration of reflection at the left boundary}
	\label{fig:reflection}
\end{figure}

The system of moment equations \eqref{eq:moments} is discretized using a first-order kinetic scheme on equidistant, structured grids (see e.g. \cite{Frank2006,Schneider2015b,Hauck2010,Garrett2014}). It is based on first discretizing  \eqref{eq:cell} in space using the first-order Upwind scheme \cite{Toro2009} and then use the method of moments in the velocity variable.
In one dimension, the scheme has the following form, abusing $v \geq 0$ for $v \in V_+$ and $v \leq 0$ for $v \in V_{-}$:
\begin{align*}
\rho_+^{n+1,i} &= \rho_+^{n,i} - \frac{\Delta t}{\Delta x} \left(q_+^{n,i}-q_+^{n,i-1}\right) + \Delta t \left(-\lambda \rho_+^{n,i} + \frac{1}{2} \lambda \rho^{n,i} + \frac{1}{4} \alpha \rho^{n,i} \Phi(\partial_x m)^{n,i}\right), \\
\rho_-^{n+1,i} &= \rho_-^{n,i} - \frac{\Delta t}{\Delta x} \left(q_+^{n,i+1}-q_+^{n,i}\right) + \Delta t \left(-\lambda \rho_-^{n,i} + \frac{1}{2} \lambda \rho^{n,i} - \frac{1}{4} \alpha \rho^{n,i} \Phi(\partial_x m)^{n,i}\right),
\end{align*}
and
\begin{align*}
q_+^{n+1,i} &= q_+^{n,i} - \frac{\Delta t}{\Delta x} \left(r_+^{n,i}-r_+^{n,i-1}\right) + \Delta t \left(-\lambda q_+^{n,i} + \frac{1}{4} \lambda \rho^{n,i} + \frac{1}{6} \alpha \rho^{n,i} \Phi(\partial_x m)^{n,i}\right), \\
q_-^{n+1,i} &= q_-^{n,i} - \frac{\Delta t}{\Delta x} \left(r_+^{n,i+1}-r_+^{n,i}\right) + \Delta t \left(-\lambda q_-^{n,i} - \frac{1}{4} \lambda \rho^{n,i} + \frac{1}{6} \alpha \rho^{n,i} \Phi(\partial_x m)^{n,i}\right),
\end{align*}
where $\rho_\pm^{0,i}$ and $q_\pm^{0,i}$ are the cell-averages of the initial conditions $\rho_\pm(0,x)$ and $q_\pm(0,x)$ in the $i$-th cell, respectively.
The chemoattractant equation \eqref{eq:chem} is discretized using an implicit finite-difference approximation, yielding
\begin{equation*}
m^{n+1,i}-D_m\frac{\Delta t}{\Delta x^2} (m^{n+1,i+1}-2m^{n+1,i}+m^{n+1,i-1}) = \beta\rho + \left(1-\delta\right) m^{n,i}.
\end{equation*}
The occurring (large and sparse) linear system is solved using the stabilized bi-conjugate gradient method implemented in the Matlab \cite{MATLAB:2012} function \verb|bicgstab|.
Analogously, the scheme extends to two dimensions, performing the ``upwinding'' on the corresponding quarter moments. Details on this implementation can be found in \cite{Ritter2015}.
We choose $\Delta t = \frac{0.5}{\frac{1}{\Delta x}+ \lambda + \alpha(s+1)}$ in one dimension and $\Delta t=\frac{0.5}{\frac{1}{\min(\Delta x,\Delta y)}+ \lambda + \alpha(s+1)}$ in two dimensions to ensure the stability of the scheme.
Since it is essential to ensure the realizability conditions in case of the minimum-entropy closure (otherwise, the defining minimization problem has no solution), we use a projector into the set of realizable moments, as derived above. We set after each new calculation of $\rho_\pm$ and $q_\pm$:
\begin{align*}
\rho_\pm &=\max(\rho_\pm, 10^{-14}), \\
q_\pm &=\pm \max(\pm q_\pm, 0), \\
q_\pm &=\begin{cases} q_\pm &\text{ if } \|q_\pm\| \leq 1, \\ \frac{q_\pm}{\|q_\pm\|} &\text{ if }  \|q_\pm\|> 1, \end{cases}
\end{align*}
and analogously in two dimensions.
Since we want to model the chemotaxis of cells in an isolated domain, we implement reflective boundary conditions for equation~\eqref{eq:cell}. Therefore we construct a layer of ghost cells, which are to be the neighbours of the boundary cells of the domain, and define their values before each step in the algorithm.
In the one dimensional case we observe the following: Only those cells that are ``travelling to the left'' are reflected from the left boundary and afterwards ``travel to the right''. That means only those cells with $v \in V_-$ are reflected on the left boundary and afterwards have $v \in V_+$. Conversely, only those cells with $v \in V_+$ are reflected from the right boundary to then have $v \in V_-$. Therefore we set:
\begin{align*}
\rho_+(\text{left ghost cells}) &= \rho_-(\text{left boundary cells}), \\
\rho_-(\text{right ghost cells}) &= \rho_+(\text{right boundary cells}),
\end{align*}
and
\begin{align*}
q_+(\text{left ghost cells}) &= -q_-(\text{left boundary cells}), \\
q_-(\text{right ghost cells}) &= -q_+(\text{right boundary cells}) \ .
\end{align*}
In two dimensions the idea remains the same, but we have two directions of reflection, see Figure~\ref{fig:reflection}.
In all cases, the chemoattractant gets reflected from the boundary in an isotropic manner. Therefore $m(Ghostcells)=m(Boundary)$. 
For Equation~\eqref{eq:chem} we use Neumann boundary conditions, which can be implemented easily by defining the neighbours of the boundary to be themselves, e.g the left neighbours of the left boundary cells to be the left boundary cells.
\section{Numerical experiments}
To test the implemented scheme, we will submit it to a number of examples and compare the results to those computed with other means. We will also interpret the results in terms of our physical model of chemotaxis. All codes are implemented in Mathworks Matlab \cite{MATLAB:2012}.

\subsection{One dimension}
In one dimension, we compare the results of our simulation to those computed by solving the kinetic equations~\eqref{eq:cell} and~\eqref{eq:chem} directly. For the latter, we used an implementation by Andreas Roth \footnote{TU Kaiserslautern, roth@mathematik.uni-kl.de}, which also uses the first-order Upwind scheme, but is explicit in the chemoattractant equation.
Unless otherwise noted we will use isotropic initial conditions, which means $f(0,x,v)=f(0,x)$ and therefore 
\begin{equation*}
\rho_\pm (0,x) = \int_{V_\pm}  f(0,x,v) dv = f(0,x) \ .
\end{equation*}
This implies
\begin{equation*}
q_\pm (0,x) = \int_{V_\pm} v f(0,x,v) dv = \pm \frac{1}{2} f(0,x) = \pm \frac{1}{2} \rho_\pm(0,x).
\end{equation*}

\subsubsection{One Spike}
\begin{figure}
	\externaltikz{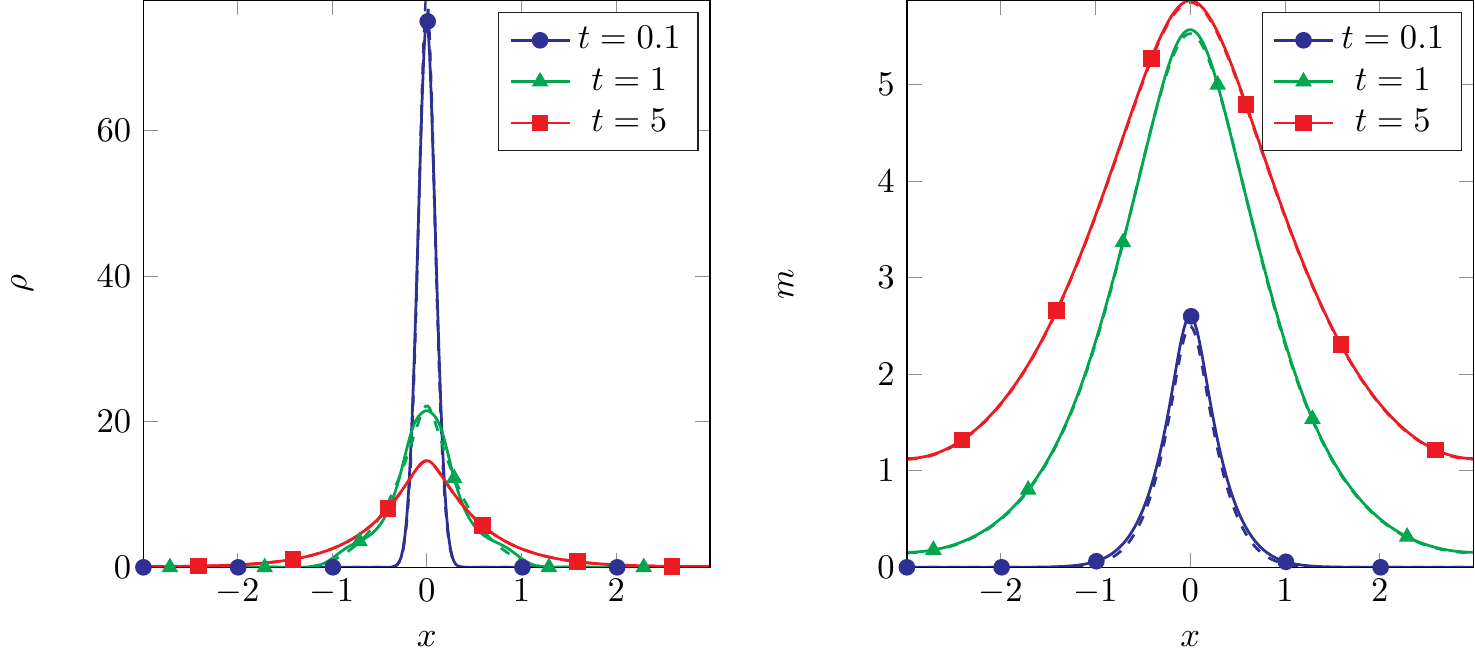}{\input{images/1D_OneSpike}}
	\caption{Half-moment minimum-entropy (solid with marks) and kinetic model (dashed) for the one-spike example.}
	\label{fig:1D_OneSpike}
\end{figure}

First, we simulate the effects of chemotaxis on an aggregation of cells in the centre of the domain in the absence of a chemoattractant initially. This is modelled by the initial data
\begin{align*}
\rho_\pm (0,x) &= \frac{1}{2} \left(100 \exp\left(-\frac{x^2}{0.01}\right) + 10^{-4}\right) \\
m (0,x) &= 0
\end{align*}
on the domain $x \in \left[-3,3\right]$ with the parameters $\alpha=2$, $D_m=1$, $\beta=1$, $\delta=1$, $\lambda=2$ and $s=0$.\\
Figure~\ref{fig:1D_OneSpike} shows, that the aggregation of $\rho$ diffuses with decreasing speed, resulting in a ``smeared out'' version of the initial state. The concentration of the chemoattractant $m$ first increases drastically until it matches $\rho$, which is caused by the production of the chemoattractant by the cells themselves, and then flattens out simultaneously with $\rho$.
Varying the parameters confirms their physical interpretations supplied in the beginning: increasing the diffusivity $D_m$, $\lambda$ speeds up the diffusion process and results in a more flattened out distribution, while increasing the parameters $\alpha$ or $\beta$ leads to reciprocal reinforcing of the spiky distribution and therefore counteracts the diffusion. \\
When comparing the results obtained by using the different closure relations, one can only observe negligible differences, since the initial condition does not provoke negativity of $\rho$ for the linear or fractional closures. Therefore we only consider the exponential closure in the following. \\
Most importantly, we see that the macroscopic and kinetic solutions behave very similarly. The minimum-entropy solution is slightly more diffusive, but with advancing time the deviation is negligible and almost invisible to the naked eye. This shows that the macroscopic model yields very good results for this example while being much cheaper computationally: the kinetic reference solution has computation times \footnote{Dual Core 2.6 GHZ, 8 GB RAM} of 9.9s and 126.8s for  $\Delta x = 0.1$ and $\Delta x=0.02$ respectively with $T=5$, whereas the macroscopic half-moment solution needs 1.0s and 3.3s with the exponential and 0.3s and 1.3s with the linear closure.

\subsubsection{Two Spikes}
In this example the (non-isotropic) initial data describes two spikes, that are located symmetrically around the centre of the domain and are moving towards it. To simplify the simulation, we consider a constant chemoattractant concentration in the shape of a downwards opening parabola, which is fixed by setting the parameters within Equation~\eqref{eq:chem} to zero. Therefore we have:
\begin{align*}
\rho_+ (0,x) &= 100 \exp\left(-\frac{(x+1)^2}{0.01}\right) + 10^{-4} \\
\rho_- (0,x) &= 100 \exp\left(-\frac{(x-1)^2}{0.01}\right) + 10^{-4}\\
q_+ (0,x) &= 100 \exp\left(-\frac{(x+1)^2}{0.01}\right) \\
q_- (0,x) &= 100 \exp\left(-\frac{(x-1)^2}{0.01}\right) \\
m (0,x) &= -x^2+9
\end{align*}
on the domain $\left[-3,3\right]$ with the parameters $\alpha=\frac{1}{2}$, $D_m=0$, $\beta=0$, $\delta=0$, $\lambda=\frac{1}{2}$ and $s=0$.  \\
For the kinetic method, we use the initial data
\begin{align*}
f(0,x,v) &= \frac{100}{0.05 \sqrt{\pi}} \exp\left(-\frac{(v+1)^2}{0.01}\right) \exp\left(-\frac{(x-1)^2}{0.01}\right) \\
&+ \frac{100}{0.05 \sqrt{\pi}} \exp\left(-\frac{(v-1)^2}{0.01}\right) \exp\left(-\frac{(x+1)^2}{0.01}\right).
\end{align*}
In this example we want to compare the half-moment solution with the full-moment $\MN[1]$ model whose initial condition can be obtained by setting $\rho = \rho_++\rho_-$ and $q = q_++q_-$.

\begin{figure}
	\externaltikz{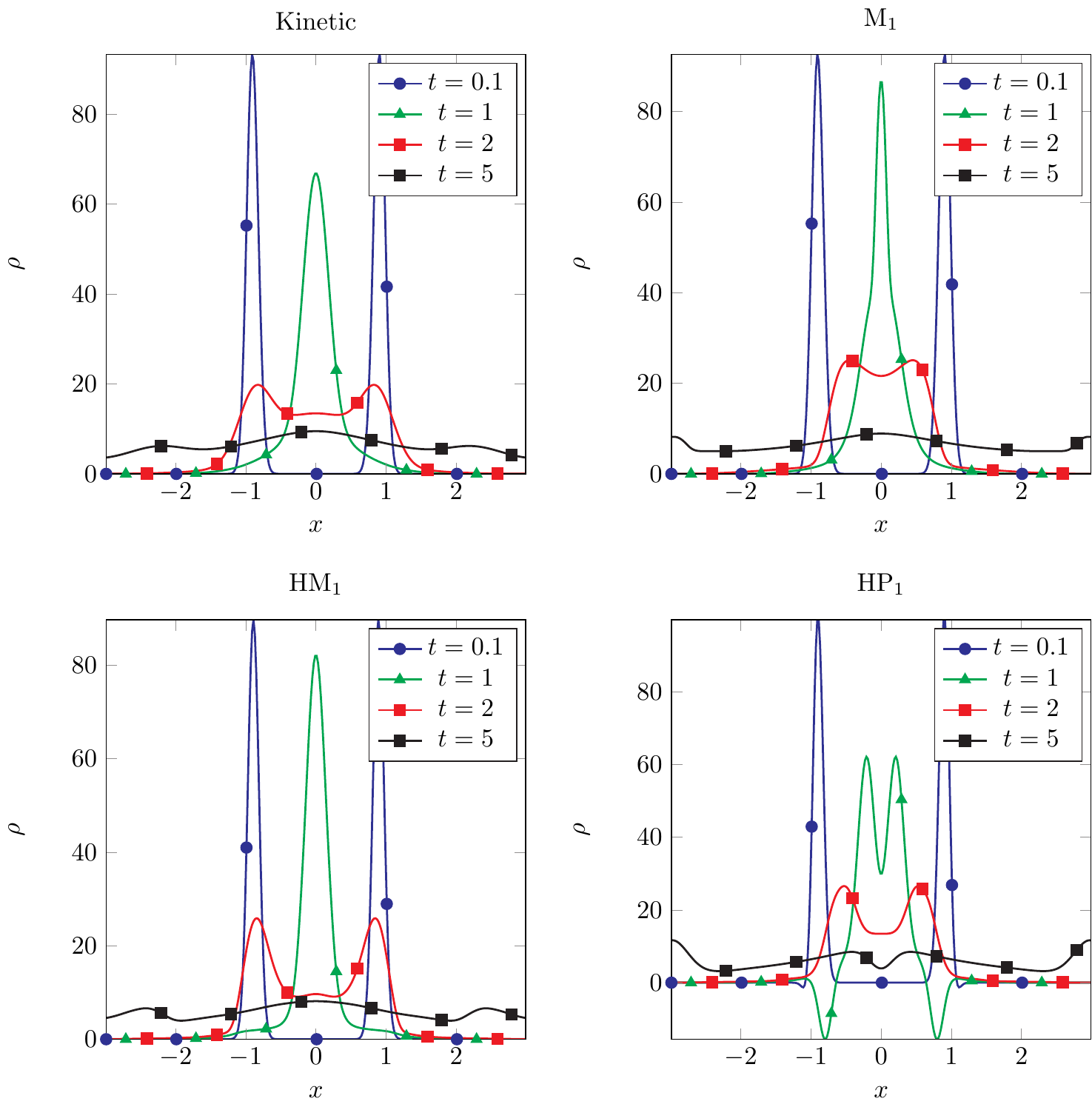}{\input{images/1D_TwoSpikes}}
	\caption{Comparison of models for the two-spikes example.}
	\label{fig:1D_TwoSpikes}
\end{figure}

We observe in Figure~\ref{fig:1D_TwoSpikes} that the two spikes move towards, collide in and then oscillate once around the centre of the domain until they reach a steady state forced by the constant chemoattractant concentration.
In this case, there are noticeable differences in the closures. In the linear half-moment $\HPN[1]$ model the density $\rho$ takes negative values on the outer edge of the spikes, but these effects cancel out again after the collision in the centre of the domain.
The full-moment minimum-entropy $\MN[1]$ model wrongly predicts a lower speed of propagation when the two spikes hit each other, caused by the well-known ``zero-netflux'' problem of this model \cite{Hauck2010,Schneider2014}.
In contrast, the speed of propagation predicted by the half-moment minimum-entropy $\HMN[1]$ model is almost exact.

\subsubsection{Kurganov Example}
\begin{figure}
	\externaltikz{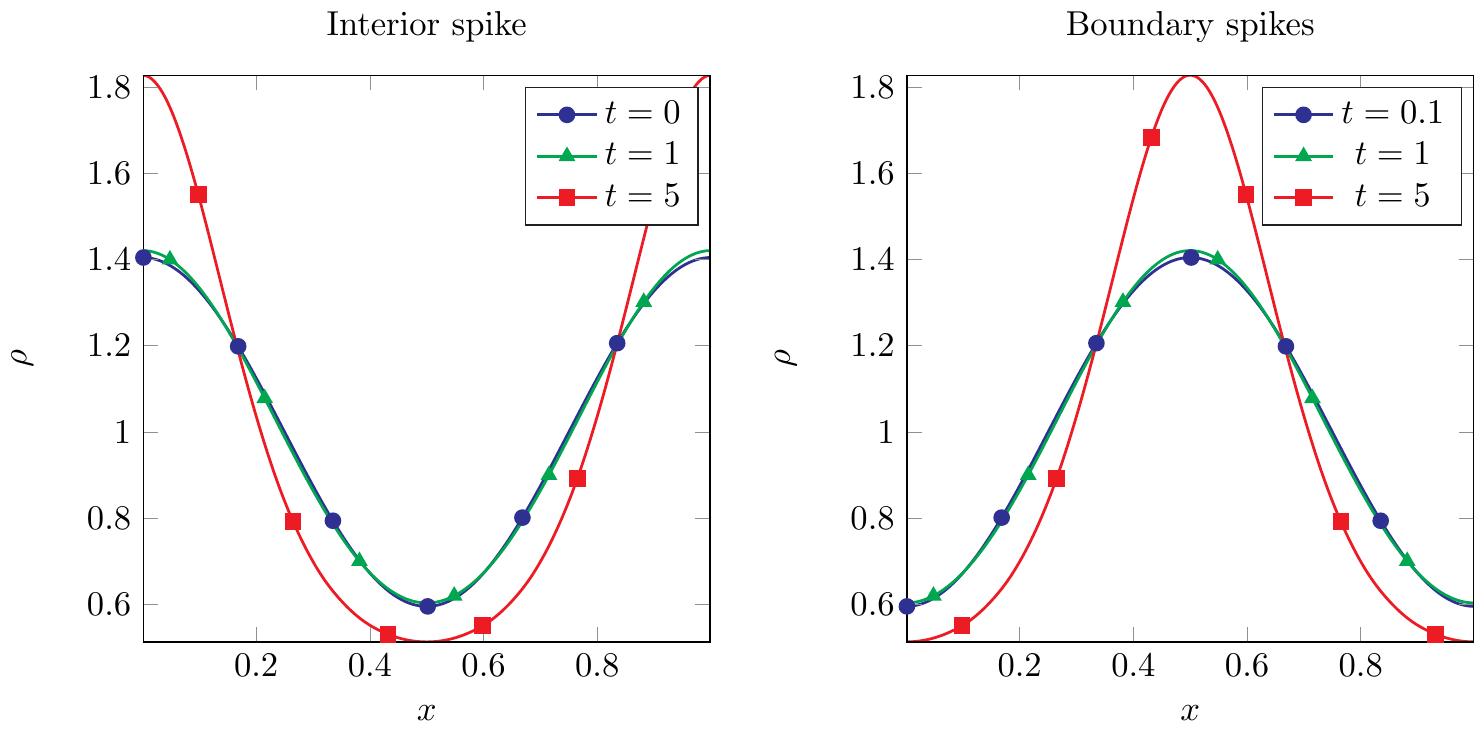}{\input{images/1D_Kurganov99}}
	\caption{Half-moment minimum-entropy model for the Kurganov examples.}
\label{fig:1D_Kurganov99}
\end{figure}

In this section we apply our simulation to an example described in~\cite{chertock} that produce bounded spiky steady states. This is of interest because it offers a model for cell aggregation which is a driving force behind embryological development and tumour formation. \\
First, we use the initial data:
\begin{align*}
\rho_\pm (0,x) &= \frac{1}{2}(1 - 0.01 (1 + 4\pi^2) \cos(2\pi x)) \\
m (0,x) &= 1 - 0.01 \cos(2\pi x)
\end{align*}
on the domain $\left[0,1\right]$ with the parameters $\alpha =1.2 (1+4\pi^2)$, $D_m=1$, $\beta=1$; $\delta=1$, $\lambda=\frac{1}{2}$ and $s=0$.
This produces an interior spike, while changing the sign to 
\begin{align*}
\rho_\pm (0,x) &= \frac{1}{2}(1 + 0.01 (1 + 4\pi^2) \cos(2\pi x)) \\
m (0,x) &= 1 + 0.01 \cos(2\pi x)
\end{align*}
yields two boundary spikes. Since the results pictured in Figure~\ref{fig:1D_Kurganov99} match those described in the paper~\cite{chertock}, we omit the comparison to the results computed with the kinetic model here. We want to remark that in this case the $\MN[1]$ model produces similar results.

\subsection{Two dimensions}
Finally, we transfer the One Spike and Two Spikes used in the one-dimensional case to two dimensions, which increases the complexity and highlights some advantages and disadvantages of the quarter-moment method. Unlike in one dimension, we do not have the means to compute a kinetic reference solution, but compare the results to those obtained with the full moment $\MN[1]$ model \cite{Klar2014,Levermore1984}. This model is well-known to be a good approximation to the true kinetic solution in many cases. All results are obtained on a grid with $\Delta x = \Delta y = 0.1$.
By default, we will use isotropic initial conditions:
\begin{equation*}
\rho_{\pm \pm'} (0,\x) = \int_{V_{\pm \pm'}}  f(0,\x,v) d\phi dr = \pi f(0,\x).
\end{equation*}

\subsubsection{One Spike}
\begin{figure}[htbp]
	\externaltikz{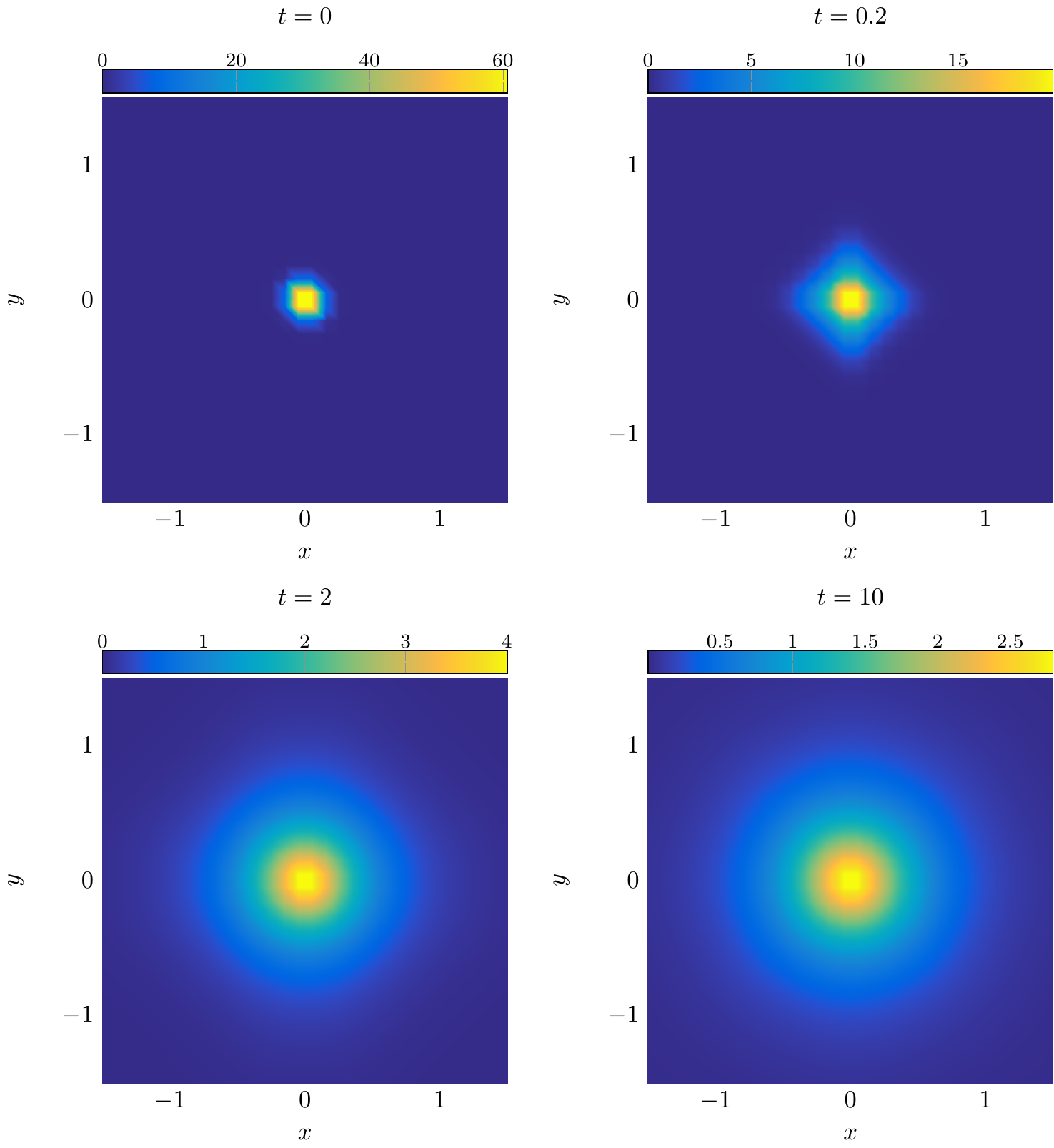}{\input{images/2D_OneSpike-60-10-expo_rho}}
	\caption{Quarter-moment $\QMN[1]$ density $\rho$ for the One Spike example.}
	\label{fig:2D_OneSpike}	
\end{figure}

\begin{figure}[htbp]
		\externaltikz{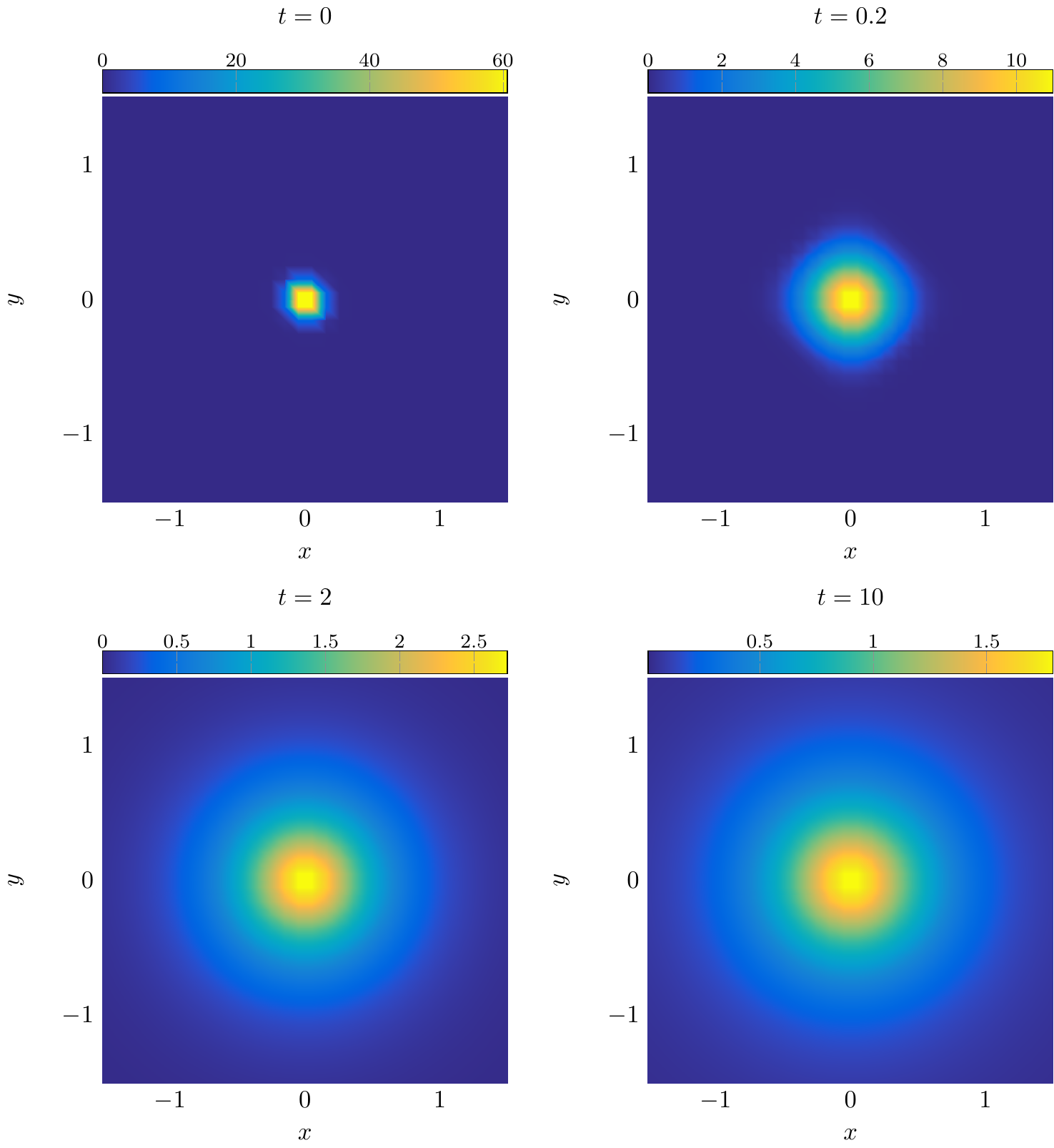}{\input{images/2D_OneSpike-60-10-expo_rho_comp}}
	\caption{Full-moment $\MN[1]$ reference density $\rho$ for the One Spike example.}		
	\label{fig:2D_OneSpike_comp}
\end{figure}

Analogously to our first example in one dimension, we start by considering the simplest case of initial data: a single spike centred in the domain with no chemoattractant present. This is modelled by 
\begin{align*}
\rho_{\pm \pm}(0,\x) &= \frac{1}{4}(100 \exp(-100 (x^2+y^2) + 10^{-4}) \\
m(0,\x) &= 0
\end{align*}
with $\x=(x,y)^T \in \left[-3,3\right] \times \left[-3,3\right]$ and the parameters $\alpha=4$, $D_m=1$, $\delta=1$, $\beta=8$ ,$\lambda=2$ and $s=0$. \\
Again, we observe in Figure~\ref{fig:2D_OneSpike} that the spike diffuses with decreasing speed.
Just like in one dimension, switching between the exponential and linear closure relation has negligible effects. \\
The fact that the reference solution computed with the full-moment method pictured in Figure~\ref{fig:2D_OneSpike_comp} is more flattened out is not surprising, since it uses the Lax-Friedrich scheme, which introduces more numerical dissipation. We also observe, that while the full-moment method conserves the rotational symmetry of the initial data, the quarter-moment method favours the axial directions, resulting in a more square-shaped solution. However, this effect vanishes for advancing times. Overall, this example shows that also in two dimensions, our simulation yields reliable results.

\subsubsection{Two Spikes}
\begin{figure}[htbp]
	\externaltikz{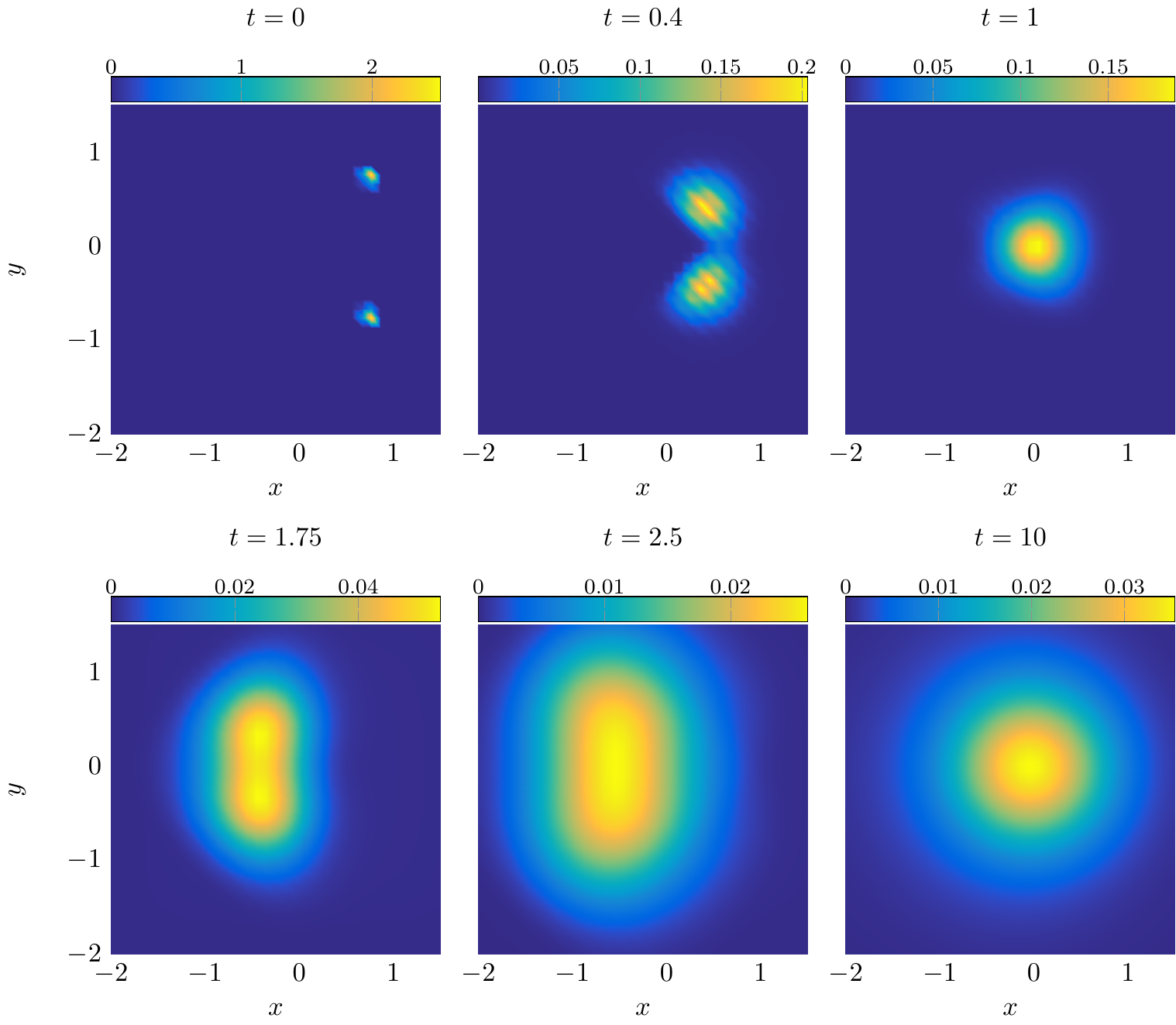}{\input{images/2D_TwoSpikesdiagonal-60-10-expo_rho_comp_seperate}}
	\caption{$\MN[1]$ superposition of the two-spikes example with diagonal placement.}		
	\label{fig:2D_TwoSpikesdiagonal_comp_seperate}	
\end{figure}

Transferring the two-spikes example to two dimensions reveals a couple of interesting effects. Again, we use the initial data of two spikes that are located symmetrically around the centre of the domain and move towards it. In two dimensions though, we can consider different alignment of the spikes, strongly influencing the outcome. 
While in one dimension the spikes opposed each other, we chose to let them propagate on orthogonal paths. In the following we distinguish the cases of diagonal and axial placement of the spikes.
To observe the effects of the simulation of Equation~\eqref{eq:cell} more closely, we again turn off all dynamics of the chemoattractant by setting the parameters in~\eqref{eq:chem} to zero.
Diagonal placement of the spikes is modelled by the initial data
\begin{align*}
\rho_{+\pm}(0,\x) &=  10^{-4} &
\rho_{- \pm}(0,\x) &= \zeta_\pm + 10^{-4} \\
q_{+ \pm}^x (0,\x) &= 0 &
q_{- \pm}^x (0,\x) &= -\frac{\zeta_\pm}{\sqrt{2}}\\
q_{+ \pm}^y (0,\x) &= 0 &
q_{- \pm}^y (0,\x) &= \pm\frac{\zeta_\pm}{\sqrt{2}}\\
m(0,\x) &= -(x^2+y^2)+18&
\zeta_\pm(0,\x) &=  100 \exp(-\frac{(x-\frac{1}{\sqrt{2}})^2+(y\pm\frac{1}{\sqrt{2}})^2}{0.001})
\end{align*}
on the domain $\left[-3,3\right] \times \left[-3,3\right]$ with the parameters $\alpha=\frac{2}{\pi}$, $D_m=0$, $\beta=0$, $\delta=0$, $\lambda=\frac{1}{\pi}$ and, $s=1$.
Similarly, the axial placement can be obtained by rotating the whole setup counter-clockwise by $45$ degrees. Note that now, the initial distribution is no longer supported in the quarter spheres, e.g.
\begin{align*}
\rho_{- -}(0,\x) &= \frac{1}{2}(100 \exp(-\frac{(x-1)^2+y^2}{0.001}) + 10^{-4}) + \frac{1}{2}(100 \exp(-\frac{x^2+(y-1)^2}{0.001}) + 10^{-4}).
\end{align*}
Since the full-moment $\MN[1]$ model is rotationally invariant, we only consider the diagonal placement. The linearity of the kinetic equation (and the physical interpretation) suggests, that the two spikes should not influence each other, so we apply the full-moment model to initial data of the two spikes simultaneously and separately, superposing the two solutions to enforce the linearity. This is similar to a pencil-beam method, see e.g. \cite{Frank07}.
Similar tests have been used in \cite{Klar2014,Schneider2015c} to investigate the $\MN[1]$ and $\QMN[1]$ model in case of fibre-laydown and radiative transfer equations.
In the following we compare the results of those four computations: 
The simulation obtained by computing both spikes separately and then superposing the results acts as our reference and is pictured in Figure~\ref{fig:2D_TwoSpikesdiagonal_comp_seperate}. 

\begin{figure}
\externaltikz{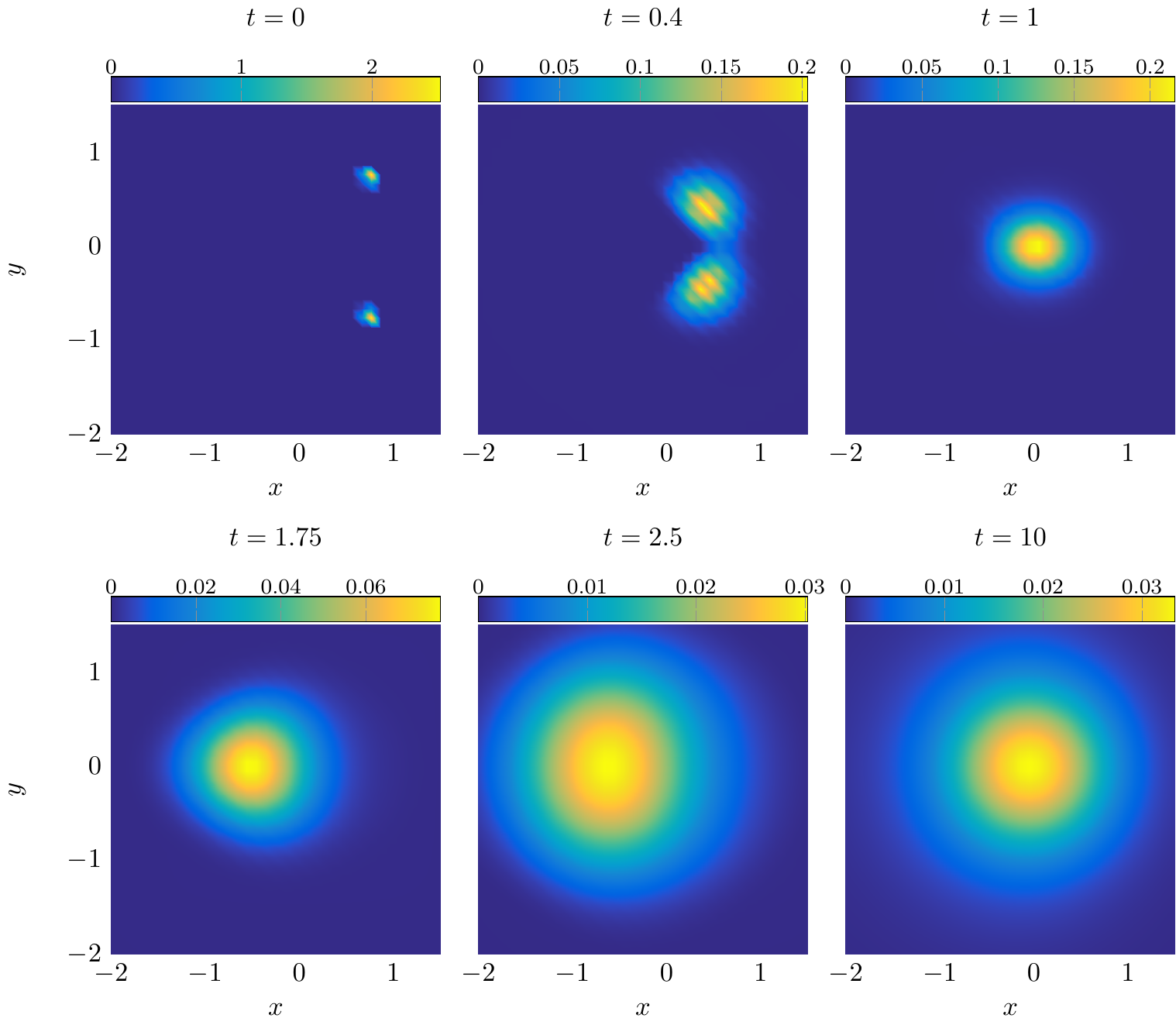}{\input{images/2D_TwoSpikesdiagonal-60-10-expo_rho_comp}}
	\caption{$\MN[1]$ solution of the two-spikes example with diagonal placement.}		\label{fig:2D_TwoSpikesdiagonal_comp}	
\end{figure}

It shows that the spikes move towards and then collide in the centre of the domain. They continue to move in their original directions but get pulled back towards the centre of the domain by the constant chemoattractant concentration, which forces a rotationally symmetric steady state.
We observe that the full moment method applied to the full initial data visualised in Figure~\ref{fig:2D_TwoSpikesdiagonal_comp}, does not produce the desired results. After colliding in the centre of the domain, the two spikes do not move on individually but rather as one before gravitating back towards the centre. 

\begin{figure}[htbp]	
\externaltikz{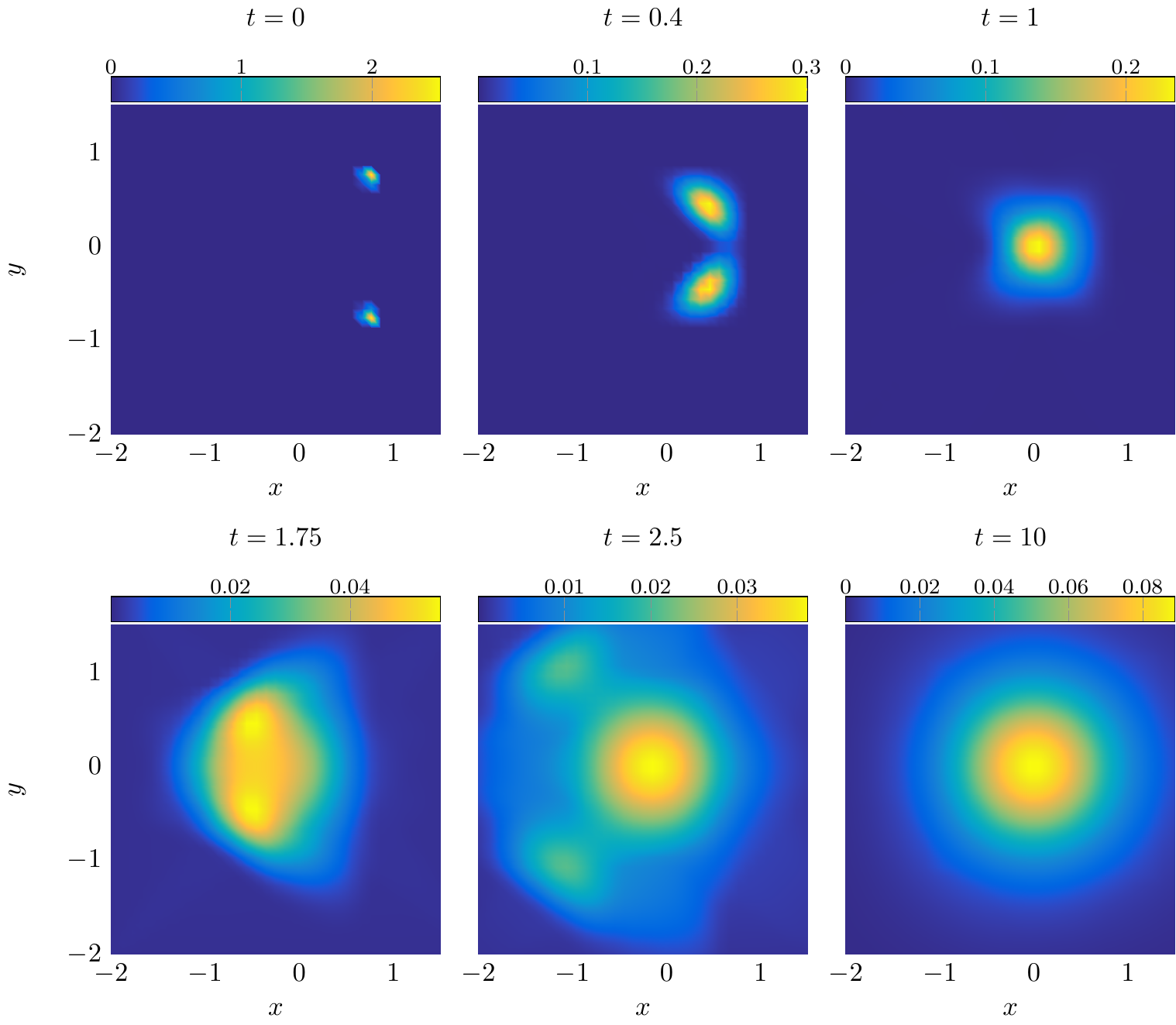}{\input{images/2D_TwoSpikesdiagonal-60-10-expo_rho}}
	\caption{$\QMN[1]$ solution of the two-spikes example with diagonal placement.}		
	\label{fig:2D_TwoSpikesdiagonal}	
\end{figure}

\begin{figure}[htbp]
\externaltikz{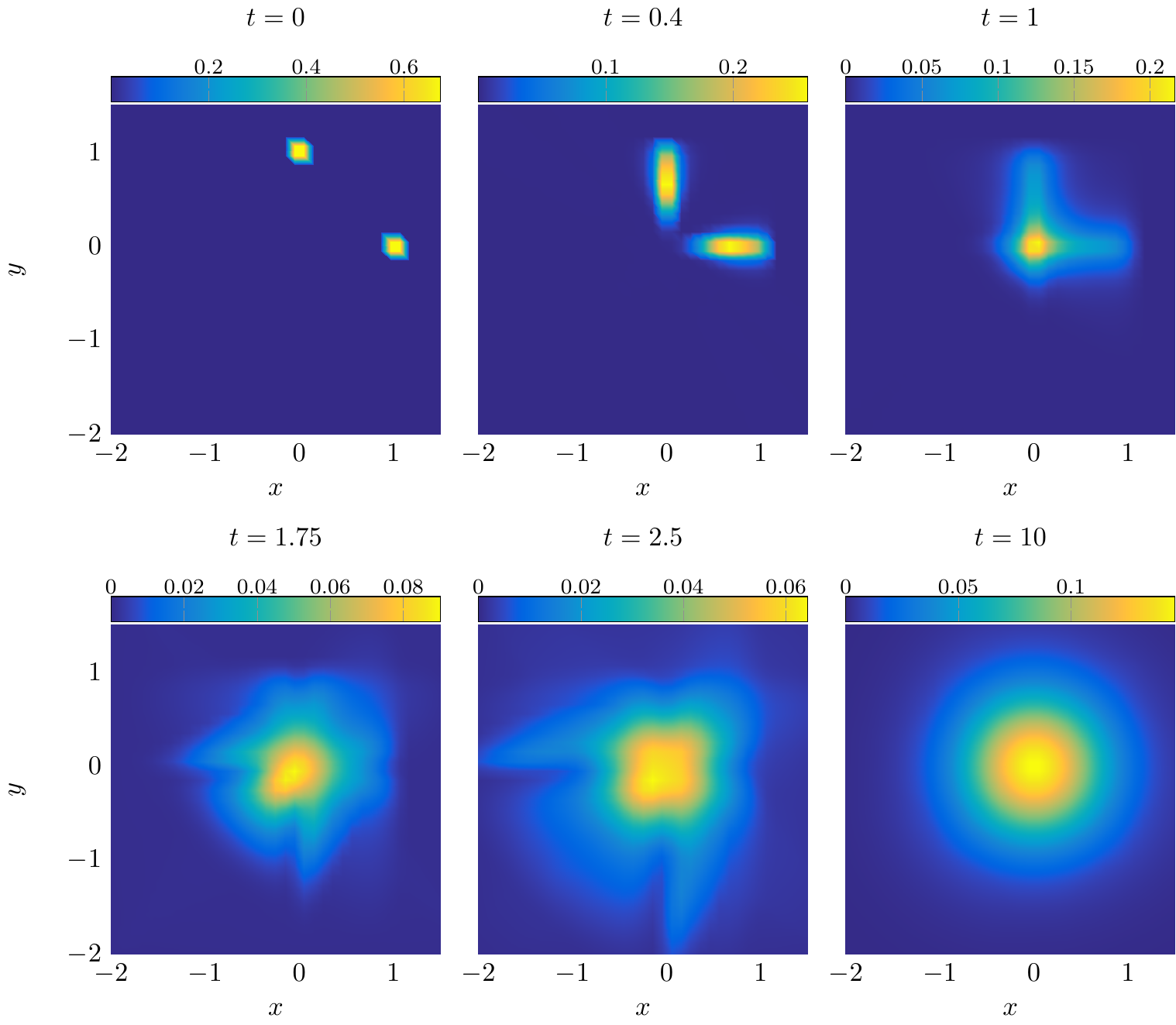}{\input{images/2D_TwoSpikesaxial-60-10-expo_rho}}
	\caption{$\QMN[1]$ solution of the two-spikes example with axial placement.}		
	\label{fig:2D_TwoSpikesaxial}		
\end{figure}

The results obtained by quarter-moment method with exponential closure are closer to the reference solution for the diagonal placement, see Figure~\ref{fig:2D_TwoSpikesdiagonal}, but not for the axial placement, see Figure~\ref{fig:2D_TwoSpikesaxial}. 

\begin{figure}[htbp]
\externaltikz{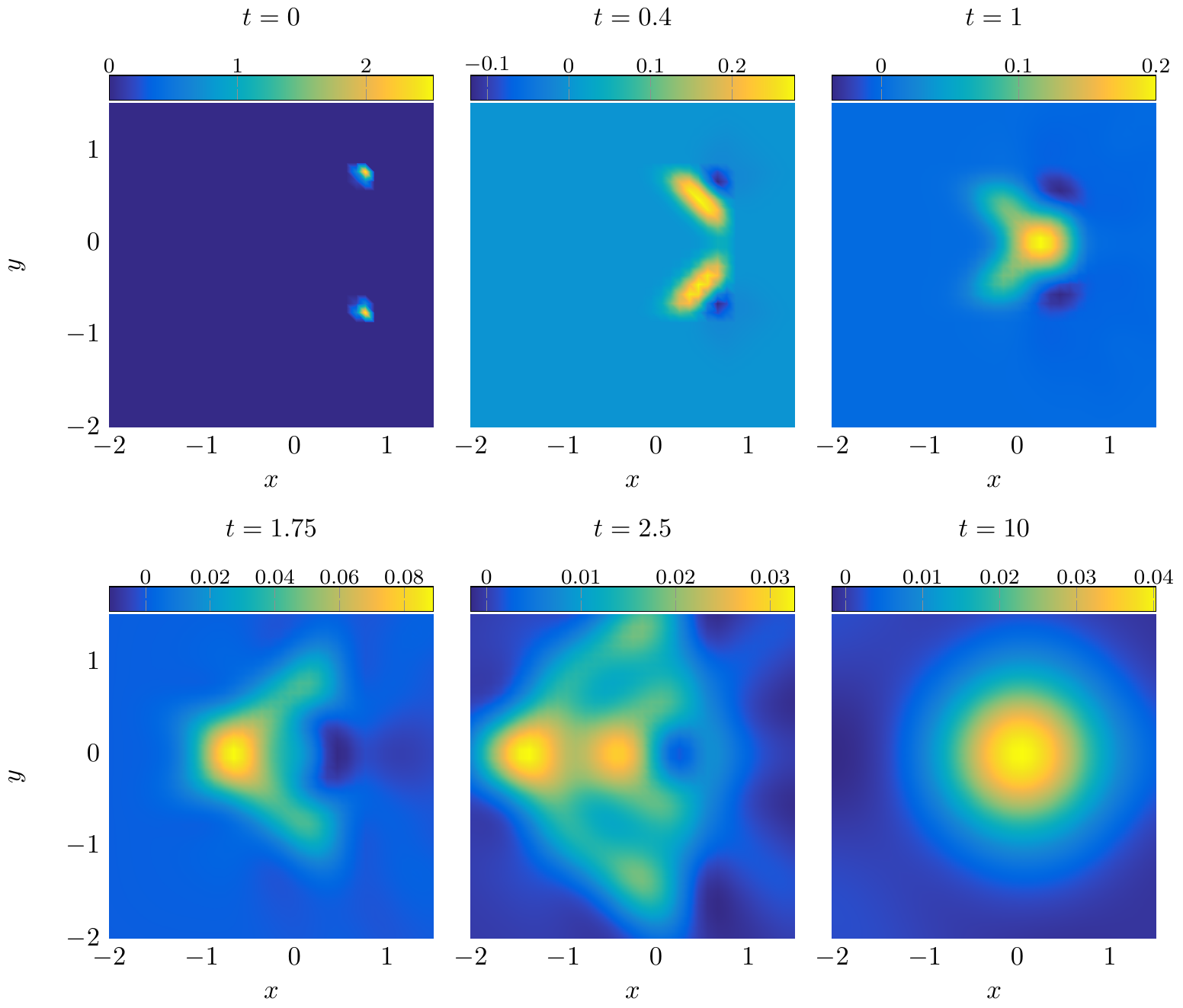}{\input{images/2D_TwoSpikesdiagonal-60-10-lin_rho}}
	\caption{Results for Two Spikes example with diagonal placement using quarter-moment method, linear closure and \mbox{$\Delta x=0.1$}}		
	\label{fig:2D_TwoSpikesdiagonal_lin}	
\end{figure}

These observations can be explained in the following way: The quarter-moment method operates on the quadrants of $V$ and can therefore distinguish between the diagonally placed spikes whose velocities lie completely in $V_{-+}$ and $V_{--}$ respectively. Therefore they do not interfere, but rather superpose as desired. This does not work for the axial placement, since in this case both spikes contribute to $\rho_{--}$. The full-moment method can not distinguish between the two spikes at all, so we observe undesirable interference. \\
This benchmark problem is also well suited to demonstrate the shortcomings of the linear $\QPN[1]$ model. It is shown in Figure~\ref{fig:2D_TwoSpikesdiagonal_lin} that the particle density $\rho$ can become negative, resulting in physically meaningless solutions.

\section{Conclusions}
We have shown in various numerical experiments that the half-/ quarter-moment method yields macroscopic models that can produce satisfying simulations of chemotaxis using simple numerical schemes for partial differential equations. We have compared the results of the simulations in one dimension to those based on the kinetic equations and found that they give a good approximation while cutting down the computational expense considerably. Comparing the simulation results in two dimensions to those obtained with a full moment method has shown that the quarter-moment method is superior when dealing with non-isotropic initial data, but has its limitations caused by the geometry of the space of admissible velocities. 

A typical strategy to solve this problem is to use higher-order models (full or partial moments). For example, the second-order $\MN[2]$ model already provides reasonably better results in case of the two-dimensional two-spikes example (see e.g. \cite{Schneider2015c} for a similar test case). However, although this model has less degrees of freedoms (six equations in total while the $\QMN[1]$ model has twelve) it is much more expensive. Due to the decoupling of the quarter-moment fluxes a two-dimensional tabulation strategy can be used. This is impossible in case of the $\MN[2]$ model, which would require a five-dimensional tabulation. Thus, Newton-like algorithms have to be applied for the $\MN[2]$ model in every space-time cell, leading to a high numerical cost, which is sometimes even more expensive than solving the kinetic equation itself. 

A possibility to further increase the accuracy of the quarter-moment ansatz is to do a further refinement of the sphere, leading to a general first-order partial-moment model, which will hopefully converge to the true kinetic solution with a small number of refinements while being numerically efficient due to the inherent structure, still enabling a cheap tabulation strategy.

\bibliographystyle{siam}
\bibliography{bib,bibliography,lit,lit2,chemotaxis_coupling}

\end{document}